\documentclass[11pt, dvips]{amsart}
\usepackage{vmargin}
\usepackage{graphicx}
\usepackage{draftcopy}
\usepackage{url}

\theoremstyle{definition}

\theoremstyle{remark}

\makeatletter
\def\romenumi{%
  \def\theenumi{\roman{enumi}}%
  \def\p@enumi{\theenumi}%
  \def\labelenumi{(\@roman\c@enumi)}}
\makeatother

\numberwithin{equation}{section}

\setpapersize{A4}
\setmarginsrb{20mm}{20mm}{20mm}{20mm}{12pt}{10mm}{12pt}{10mm}

\DeclareMathOperator{\vol}{vol}
\newcommand{\Em}[1]{\textbf{#1}}

\begin{document}

\title[] {Probability calculations under the IAC hypothesis}
\author{Mark C. Wilson}
\address{Department of Computer Science, University of Auckland,
 Private Bag 92019 Auckland, New Zealand}
\email{mcw@cs.auckland.ac.nz}
\author{Geoffrey Pritchard}
\address{Department of Statistics, University of Auckland, Private Bag 92019
 Auckland, New Zealand}
\email{geoff@stat.auckland.ac.nz}

\date{\today}



\begin{abstract}
We show how powerful algorithms recently developed for counting lattice
points and computing volumes of convex polyhedra can be used to compute
probabilities of a wide variety of events of interest in social choice
theory. Several illustrative examples are given.
\end{abstract}

\maketitle

\section{Introduction}
\label{sec:intro}

Much research has been undertaken in recent decades with the aim of
quantifying the probability of occurrence of certain types of election
outcomes for a given voting rule under fixed assumptions on the
distribution of voter preferences. Most prominent among these outcomes
of interest are the so-called voting paradoxes, which have been shown to
be unavoidable, hence the interest in how commonly they may occur. The
survey \cite{GeLe2004} discusses these questions
and gives a summary of results up to 2002.

In very many cases, particularly under the IAC hypothesis on voter
preferences, the calculations involved amount simply to counting integer
lattice points inside convex polytopes. In the social choice literature,
two main methods have been used to carry out such computations. The
first, dating back several decades, decomposes the polytope into smaller
pieces each of which can be treated by elementary methods involving
simplification of multiple sums. This method works fairly well for
simple problems but requires considerable ingenuity and perseverance to
carry out even for moderately complicated ones. More recently, more
powerful methods have been introduced in \cite{HuCh2000,
Gehr2002b} but there are several recent instances where even these
methods did not suffice to solve natural questions about 3-candidate
elections.

The purpose of the present paper is to point out that there is an
established mathematical theory of counting lattice points in convex
polytopes (and the closely related issue of computing the volume of such
a region), which has been partially rediscovered by workers in social
choice theory. The area has recently been the subject of active research
(see \cite{Delo2005} for a good summary). Several more efficient
new algorithms have been devised and implemented in publicly available
software.

We aim to apply these new methods to answer questions in voting theory
that have proven beyond the reach of previous authors.  In addition we
corroborate, correct, and unify the derivation of some previously
published results by using this methodology. We believe that the solution
of many hitherto difficult problems can now be relegated to a trivial
computation. This should open the way for social choice theorists to
tackle more difficult and realistic problems. We note that 
Lepelley, Louichi and Smaoui \cite{LLS2006} have recently, and independently from us, 
circulated a preprint with a similar goal, which covers very similar ground. The
fact that two groups of researchers discovered this approach almost
simultaneously shows that the time has indeed come for these methods to
be assimilated by the social choice community.

The basic idea is that many sets of voting situations that are of
interest can be characterized by linear equations and inequalities. The
variables are usually the numbers of voters with each of the $m!$
possible preference orders, where $m$ is the number of alternatives. The
set of such (in)equalities defines a convex polytope in $\mathbb{R}^d$
for some $d$, given by $Ax \leq b$ for some matrix $A$ (here $d \leq m!$
and the inequality may be strict, because we may first use equality
relations to eliminate variables and reduce dimension). Each lattice
point will correspond to a voting situation in the desired set. The
probability that a randomly chosen situation has the property under consideration
 is therefore a straightforward ratio of lattice point counts. Dividing through by
$n$, the total number of voters, yields a convex polytope $P$, independent of $n$,
 in $\mathbb{R}^d$. For a given number $n$ of voters, the
dilation $nP$ describes the set of lattice points that we wish to
enumerate.

\section{Counting lattice points in convex polytopes}
\label{sec:convex}
We give only a brief description here. For more information we recommend
\cite{Delo2005}. 

The \Em{Ehrhart series} of the rational polytope $P$ is a rational
generating function $F(t)= P(t)/Q(t) = \sum_n a_n t^n$ whose $n$th
Maclaurin coefficient $a_n$ gives the number of lattice points inside
the dilation $nP$. The function $f:n \mapsto a_n$ is known to be a
polynomial of degree $d$ if all the vertices of $P$ are integral;
otherwise it is a quasipolynomial of some minimal period $e$. That is,
the restriction of $f$ to each fixed congruence class modulo $e$ is a
polynomial.

It is known that $e$ is a divisor of $m$, where $m$ is an integer such
that  all coordinates of vertices of $mP$ are integers. The
least such $m$ is the least common multiple of the denominators of the
coordinates of the vertices of $P$ when each coordinate is written in
reduced terms. However there are examples where $e < m$ \cite{McWo2005}. 
A method for determining $e$ was presented in \cite{HuCh2000}.

Many questions in voting theory are of most interest in the asymptotic
case where $n \to \infty$. For small $n$, issues such as the method of
tiebreaking used assume great importance, whereas in the limit such
issues disappear (the situations in which ties occur correspond in the
limit to the boundary of $P$). We focus on limiting results in the
present paper.

The leading coefficient of the quasipolynomial $f$ is the same for all
congruence classes: only the lower degree terms differ. It is well known
that this leading coefficient is precisely the volume of $P$. For many
purposes, knowledge of this coefficient is sufficient. The limiting
probability under IAC as $n \to \infty$ is simply the volume of $P$
divided by the volume of $X$ where $X$ is the analogously defined
polytope that describes all possible voting situations.

To compute the number of lattice points in $nP$, if that amount of
detail is desired, we may use one of several algorithms.  An attractive
approach pioneered by Barvinok makes heavy use of rational generating
functions; this is implemented in the software {\tt LattE}
\cite{DHTY2004, latte}. There are also several algorithms available
for volume computation; see \cite{BEF2000} for a survey of
algorithms, a hybrid of which has been implemented in {\tt vinci}
\cite{vinci} for floating point computation only. One of these
algorithms has been used in the Maple package {\tt Convex}
\cite{convex-package}, and uses exact rational arithmetic.

The software {\tt LattE} gives the Ehrhart series as standard output. In
order to extract the quasipolynomial formula for $f(n)$ from the Ehrhart
series, we may use interpolation. On each congruence class modulo $e$,
we must evaluate $f(n)$ at $d+1$ distinct values of $n$ in this class.
Given the explicit expression $F(t) = P(t)/Q(t)$ and a computer algebra
system, such evaluations are trivially obtained (the $a_n$ satisfy a
linear recurrence relation with constant coefficients). The Lagrange
interpolation formula then yields the desired formula for the particular
polynomial that is applicable for the given congruence class.

Another (generally less efficient) method of extraction is to decompose $F(t)$ 
into partial fractions. Note that $F(0) = 1$ and we can arrange so that $Q(t)$
factors as $\prod_j (1 - \alpha_j t)$ for some complex numbers
$\alpha_j$, possibly not distinct. We then have the partial fraction decomposition 
$$F(t) = \sum_\alpha \sum_k c_{\alpha,k} (1 - \alpha t)^{-k}$$
where $\alpha$ runs over the roots of $Q$ and $k$ runs from $1$ to the multiplicity of $\alpha$.

This shows how the periodicity occurs: the
factorization of $Q(t)$ will introduce complex roots of unity and the
terms corresponding to powers of these will simplify on each congruence class. 
In fact, on extracting the coefficient of $t^n$ we obtain
$$
[t^n] F(t) =  \sum_{\alpha, k}  \alpha^n c_{\alpha,k} \binom{n+k-1}{k-1}
$$
and the terms $\alpha^n$ simplify on each equivalence class modulo $e$.
Note that $e = 1$ (that is, $f(n)$ is a single polynomial) if and only if $Q$
factors completely over the rationals.

Note that since we know \textit{a priori} that the coefficient of $t^n$ is
polynomially growing, all $\alpha$ with $|\alpha| \neq 1$ can be
ignored, since their contribution must cancel (otherwise we would obtain
terms exponentially growing or decreasing in $n$). Unfortunately this
observation does not help in the present case, because the Ehrhart
series has a denominator of the form $\prod_i (1 - t^{a_i})$,  so all the
$\alpha_j$ above are in fact roots of unity.

In summary, the Ehrhart series contains all information required to
solve the problem of counting lattice points in polytopes parametrized
by a single parameter $n$. The hardest step is usually determining the
 minimal period $e$.

\section{Examples}
\label{sec:examples}

In this section we compute, using the recipe above, a few probabilities
under IAC that have been considered in the recent social choice
literature. We emphasize problems where older methods have not yielded
an answer, but also check results obtained by previous authors using
older methods. Some of these earlier results appear to be incorrect. The 
use of a computer algebra system such as Maple \cite{maple} is essential for 
some of the more complicated examples.

\subsection{Manipulability}
\label{ss:manip}

We first consider the probability under IAC that a voting situation in a
$3$-candidate election is manipulable by some coalition. Counterthreats
are not considered --- we assume that some group of voters with incentive
to manipulate will not be opposed by the other, naive, voters. See
\cite{PrWi2005, FaLe2006} for more discussion of these
(standard) assumptions.

For the classical rules plurality and antiplurality, the answer is
known: 7/24 and 14/27 respectively \cite{LeMb1987,
LeMb1994}. These results were derived by the earliest methods
described above and required considerable hand computation. 
However, for the Borda rule, no such result has been
derived even using more sophisticated methods. A good numerical
approximation to the limit has been obtained. In
\cite{FaLe2006} the authors used the method of
\cite{HuCh2000} to obtain bounds on the solution but could not carry
out the full computation. Using their method requires interpolation,
hence computing the first $6e$ coefficients of the Ehrhart series, where
$e$ is the minimal period of the quasipolynomial. They showed that $e >
48$, and since they computed these coefficients by exhaustive
enumeration, it was not possible to carry out the computation to the end
(the number of voting situations is of order $n^5$). They estimated a
value of $0.5025$ for the limit.

However with more powerful tools the answers are easily obtained. We let $n_1, \dots
,n_6$ denote the number of voters with sincere preference order $abc,
acb, bac, bca, cab, cba$ respectively, and let $x_i = n_i/n$. Then
$\sum_i x_i = 1$ and $x_i \geq 0$.  We use the linear systems derived for 
general positional rules in \cite{PrWi2005}. As shown in \cite{PrWi2005} we may assume 
without
loss of generality that $a$ wins, $b$ is second, and $c$ last in the
election (this assumption will only affect lower order terms in our resulting quasipolynomial, 
and this is inevitable when different tie-breaking assumptions are made). 
Thus we must multiply our final answer by $6$ since we are only considering
one of the $3!$ equally likely permutations of the candidates. 

\subsubsection*{Plurality}

We first consider the plurality rule. We define polytopes $P_b, P_c, P_{bc}$ as follows. 
Consider the inequalities
\begin{align}
\label{eq:plur a-b}
0 & \leq x_1 + x_2 - x_3 - x_4   \qquad \text{($a$ beats $b$ (sincere))} \\
\label{eq:plur b-c}
0 & \leq x_3 + x_4 - x_5 - x_6   \qquad \text{($b$ beats $c$ (sincere))} \\
\label{eq:plur-str b-a}
0 & \leq -x_1 -x_2 + x_3 + x_4 + x_6 \qquad \text{($b$ beats $a$ (strategic))}\\
\label{eq:plur-str b-c}
0 & \leq -x_1 -x_2 + 2x_3 + 2x_4 -x_5 + 2x_2 \qquad \text{($b$ beats $c$ (strategic))}.
\end{align}
The polytope $P_b$ (the region where manipulation in favour of $b$ is possible) is defined by 
the inequalities
\eqref{eq:plur a-b} -- \eqref{eq:plur-str b-c}, the equality $\sum_i x_i = 1$, and the condition
 that all $x_i$ are nonnegative. Polytope $P_c$ is obtained by
applying the permutation $b \leftrightarrow c$, which induces the
permutation $x_1 \leftrightarrow x_2, x_3 \leftrightarrow x_5, x_4
\leftrightarrow x_6$, and $P_{bc} = P_b \cap P_c$ is just given by the union of the two
sets of inequalities defining $P_b$ and $P_c$.

The software {\tt LattE} readily computes the Ehrhart series of each polytope. They are

\begin{align*}
H_b & =  {\frac {12\,{t}^{12}+24\,{t}^{11}+44\,{t}^{10}+56\,{t}^{9}+66\,{t}^{8}+64\,{t}^{
7}+63\,{t}^{6}+44\,{t}^{5}+30\,{t}^{4}+14\,{t}^{3}+6\,{t}^{2}+2\,t+1}{ \left( 1
-t \right) ^{2} \left( 1-{t}^{3} \right) ^{4} \left( 1+t \right) ^{4} \left( 1+{t
}^{2} \right) ^{3}}}
\\
H_c & = {\frac {8\,{t}^{12}+16\,{t}^{11}+26\,{t}^{10}+34\,{t}^{9}+38\,{t}^{8}+40\,{t}^{7
}+41\,{t}^{6}+30\,{t}^{5}+20\,{t}^{4}+10\,{t}^{3}+4\,{t}^{2}+2\,t+1}{ \left( 1
-{t}^{4} \right) ^{3} \left( 1-t \right) ^{2} \left( 1-{t}^{2} \right)  \left( 
1+t+{t}^{2} \right) ^{4}}}
\\
H_{bc} & =
{\frac {4\,{t}^{8}+5\,{t}^{6}+4\,{t}^{5}+4\,{t}^{4}+4\,{t}^{3}+2\,{t}^{2}+1}{
 \left( 1-t \right) ^{4} \left( 1-{t}^{4} \right) ^{2} \left( 1+t+{t}^{2}
 \right) ^{4}}}
\end{align*}

The series we require is therefore
\begin{align*}
H& :=H_b + H_c - H_{bc} \\
& = 
{\frac {16\,{t}^{12}+32\,{t}^{11}+57\,{t}^{10}+68\,{t}^{9}+78\,{t}^{8}+74\,{t}^{
7}+73\,{t}^{6}+50\,{t}^{5}+33\,{t}^{4}+14\,{t}^{3}+6\,{t}^{2}+2\,t+1}{ \left( 1
-{t}^{4} \right) ^{3} \left( 1-t \right) ^{2} \left( 1-{t}^{2} \right) 
 \left( 1+t+{t}^{2} \right) ^{4}}}.
\end{align*}
Note that in order to factor the denominator of $H$ completely we require both a cube root and 
fourth root of $1$, hence a field extension of degree 12. Thus we expect the period of 
the quasipolynomial $f(n):=[t^n]H(t)$ to be 12. We may determine the polynomial formula for $f$ 
on each congruence class in more than one way, as described in section~\ref{sec:convex}. 

First, we try interpolation. Consider the polynomial expression valid for $f(n)$ when 
$n \equiv 0 \mod 12$. This is a polynomial of degree 5 in $n$. We compute the values $f(12j)$ 
for $j=0, \dots, 5$ and then determine the unique interpolating polynomial of degree $5$ 
determined by these points, via, say, the Lagrange 
inversion formula. The built-in commands in Maple find this polynomial immediately: the answer is 
$$
f(n) = {\frac {7}{17280}}\,{n}^{5}+{\frac {1}{108}}\,{n}^{4}+{\frac {3}{32}}\,{n}^{3}+{
\frac {15}{32}}\,{n}^{2}+{\frac {137}{120}}\,n+1  \qquad (n \equiv 0 \mod 12).
$$
As a check, we substitute $n = 96$ into this expression --- the correct answer, namely 
$[t^{96}] H(t) = 4176821$, is obtained. Analogous formulae can be obtained in the same way for the 
other congruence classes modulo $12$. For example, the result for $n$ congruent to $6$ modulo 
$12$ is 
$$
f(n) = {\frac {7}{17280}}\,{n}^{5}+{\frac {1}{108}}\,{n}^{4}+{\frac {3}{32}}\,{n}
^{3}+{\frac {15}{32}}\,{n}^{2}+{\frac {61}{60}}\,n+5/8 \qquad (n \equiv 6 \mod 12),
$$
while that for $n$ congruent to $1$ is given by 
$$
f(n) = {\frac {7}{17280}}\,{n}^{5}+{\frac {1}{108}}\,{n}^{4}+{\frac {341}{5184}}
\,{n}^{3}+{\frac {5}{36}}\,{n}^{2}-{\frac {917}{17280}}\,n-{\frac {209}{
1296}} \qquad (n \equiv 1 \mod 12).
$$

Note that since the number of voting situations is given by 
$\binom{n+5}{5} = (n+1) \cdots (n+5)/120$, and we have only counted
one-sixth of the manipulable situations, the limiting probability of manipulability is 720 
times the leading coefficient of $f$, namely $7/24$. This agrees with the results obtained in 
\cite{LeMb1987}. Note that the expressions for finite $n$ do not agree, probably because 
of different tie-breaking assumptions yielding slightly different sets of manipulable voting 
situations. We use random tiebreaking as described in \cite{PrWi2005}, with a winner being 
chosen uniformly at random from the 
set of those with highest score; the alternative used in many papers 
breaks the symmetry by breaking ties in favour of a fixed but arbitrary order on the candidates. 
It is clear from the discussion at the beginning of the 
proof of \cite[Theorem 2]{LeMb1987} that the latter tiebreaking method is used in that paper.

As mentioned in section~\ref{sec:convex}, another method would be to compute the full partial 
fraction decomposition of $H$ over the extension field of $\mathbb{Q}$ 
generated by a primitive 12th root of $1$. This can 
be done easily by Maple. However the result is somewhat messy and the ensuing
computation involving binomial coefficients is certainly 
no easier than using interpolation, so we omit it.

\subsubsection*{Borda}

We now consider the Borda rule. We can attempt an analysis similar to the above (the polytopes 
are defined in a similar manner, and all coefficients lie in  $\{0, \pm 1, \pm 2, \pm 3\}$), 
but we run into serious complexity issues in this case.
 
The Ehrhart series $F_b, F_c, F_{bc}$
given by {\tt LattE} are such that when $F:=F_b + F_c - F_{bc}$ is
simplified, its denominator is a product of cyclotomic polynomials
(minimal polynomials for roots of unity). The corresponding roots of
unity required are of orders whose least common multiple is $2520$.
So we are still faced with the major task of computing $e$. It is still
an open problem as to whether there exists an algorithm to determine $e$ 
which runs in polynomial time (in the input size) when the dimension is fixed. A
polynomial time algorithm to determine whether an integer $p$ is equal
to $e$ was presented in \cite{Wood2005}, but has not been
implemented in software as far as we are aware. 
Of course, we do not need to know the exact value of $e$, and we could 
assume it to be $2520$. In order to determine exact formulae for $f(n)$ in all cases by 
interpolation, we would require the first $15120$ values of $f(n)$. Trying this in Maple we obtain 
an overflow error. However it would be possible in principle to compute these using the 
recurrence supplied by the rational form of $F$. We do not proceed further along these lines, but 
we indicate how the computation would go. Writing 
$P(t) = \sum_k b_k t^k, Q(t) = \sum_k c_k t^k, F(t) = P(t)/Q(t) = \sum_n a_n t^n$ 
and comparing coefficients, we obtain $b_n = \sum_{0\leq k\leq n} c_k f(n-k)$. This constant 
coefficient linear recurrence allows us to 
determine sequentially $f(0), \dots, f(r)$ where $r = \deg P$, and for $n > r$ we have the 
defining recurrence $\sum_{0\leq k\leq n} c_k f(n-k) = 0$. In the present case 
$\deg P = 75$ and  $\deg Q = 82$, so the computation would be rather involved.

However, we can certainly determine the leading term of the quasipolynomial $f$, namely 
the volume of a certain region. It is convenient to 
eliminate $x_6$ throughout, using the sole equality 
constraint $\sum_i x_i = 1$. In other words we look at the projection onto the
subspace $x_6 = 0$. Since we are dividing by the volume of the
projection of the simplex the exact scale factor is unimportant. This
projection is defined by the conditions $x_i \geq 0$ and $\sum_{i=1}^5 x_i
\leq 1$ --- we call these the \Em{standard inequalities}. The volume in
$\mathbb{R}^5$ of this simplex is easily computed to be $1/5! = 1/120$. Recalling the 
factor of $6$ mentioned above, we shall therefore multiply the volume 
answer obtained below by $720$ to compute the limiting probability.

The volume required is given by inclusion-exclusion as $\vol(R_b) +
\vol(R_c) - \vol(R_b \cap R_c)$ where $R_b, R_c$ respectively denote the
region for which manipulation in favour of $b$ or $c$ is possible. 

The conditions describing the sincere outcome reduce, after elimination
of $x_6$, to 
\begin{align}
\label{eq:borda a-b}
2x_1 + 3x_3 + 2x_4 - x_5 & \geq 1;  \qquad \text{($a$ beats $b$ (sincere))} \\
\label{eq:borda b-c}
2x_1 + 3x_2 - x_4 + 2x_5 & \geq 1;  \qquad \text{($b$ beats $c$ (sincere))} 
\end{align}
while the conditions describing the outcome after manipulation amount to
\begin{align}
\label{eq:borda-str b-a}
3x_1 + 4x_2 + 3x_5 & \leq 2  \qquad \text{($b$ beats $a$ (strategic))}\\
\label{eq:borda-str b-c}
x_1 + 2x_2 + 2x_5 & \leq 1 \qquad \text{($b$ beats $c$ (strategic))}.
\end{align}
Now $R_b$ is defined by the standard inequalities and those in
\eqref{eq:borda a-b} -- \eqref{eq:borda-str b-c}. Also $R_c$ is obtained by
applying the permutation $b \leftrightarrow c$, which induces the
permutation $x_1 \leftrightarrow x_2, x_3 \leftrightarrow x_5, x_4
\leftrightarrow x_6$, and $R_{bc}$ is given by the union of the two
sets of inequalities defining $R_b$ and $R_c$.

The package \texttt{Convex} \cite{convex-package} immediately yields the
answer when given this input. The respective volumes of $R_b, R_c,
R_{bc}$ are $371/559872, 881/6531840, 170873/1714608000$ and the
required limit is precisely $132953/264600 \approx 0.5024678760$.

The large denominators in the fractions above give a clue to the
difficulty of this problem.  \texttt{Convex} also computes the
vertices of the polytope. The least common multiple of the denominators
of the coordinates of the vertices is $72$ for $R_b$, $504$ for $R_c$,
$1260$ for $R_{bc}$. Thus the minimum period $e$ is a divisor of $2^3
\cdot 3^2 \cdot 5 \cdot 7 = 2520$, as we already knew from above. 

\subsection{Condorcet phenomena}
\label{ss:cond}

See the two surveys and recent book by Gehrlein \cite{Gehr1997,
Gehr2002a, Gehr2006} for more information about previous work on this
topic.

In \cite{GeLe2004} Gehrlein and Lepelley state ``A
very large number of studies (probably more than $50\%$ of the studies
that have been devoted to probability calculations in social choice
theory) have been conducted to develop representations for the
probability that Condorcet's Paradox will occur, and for the Condorcet
efficiency of various rules, with the assumptions of IC and IAC."

\subsubsection*{Condorcet's paradox}

\Em{Condorcet's Paradox} occurs in a voting situation when there is no
Condorcet winner --- that is, no one candidate beats all others when
only pairwise comparisons are considered. This occurrence is independent
of the voting rule being used. To compute its likelihood, we compute the
complementary event.

Suppose that we have $3$ alternatives $a, b, c$. Let $C$ be the event
that $a$ is the Condorcet winner. This yields inequalities that boil
down to 
\begin{align}
\label{eq:cond a-b}
2x_1 + 2x_2 + 2x_3 & \geq 1  \qquad\text{($a$ beats $b$ pairwise);} \\
\label{eq:cond a-c}
2x_1 + 2x_2 + 2x_5 & \geq 1  \qquad\text{($a$ beats $c$ pairwise).}
\end{align}
Let $P_C$ be the polytope defined by these and the standard
inequalities. Then {\tt Convex} yields $\vol(P_C) = 1/384$, so that
Condorcet's Paradox occurs with asymptotic probability $1 - 3 \cdot
5!/384 = 1/16$ for IAC with 3 alternatives. 
This is of course a known result dating back several decades.

\subsubsection*{Condorcet efficiency}

Similarly we may compute the \Em{Condorcet efficiency} of a given  rule,
namely the conditional probability that it elects the Condorcet winner
given that this winner exists. For a given scoring rule defined by
weights $(1, \lambda, 0)$, let $X_\lambda$ be the event that $a$ is the
winner when this rule is used. Clearly $\Pr(X_\lambda) = 1/3$.

These conditions describing $X_\lambda$ amount to 
\begin{align}
\label{eq:gen a-b}
x_1 + (1 + \lambda) x_2 + (2 \lambda - 1) x_3 + (\lambda - 1) x_4 + 2\lambda x_5 & \geq \lambda 
\qquad \text{($a$ beats $b$ with rule $\lambda$)}  \\
\label{eq:gen a-c}
2x_1 + (2 - \lambda) x_2 + (1 + \lambda) x_3 + (1 - \lambda) x_4 + \lambda x_5 & \geq 1 
\qquad \text{($a$ beats $c$ with rule $\lambda$)}.
\end{align}

The Condorcet efficiency of rule $\lambda$ is $\Pr(X_\lambda \cap
C)/\Pr(C)$ which equals $3 \cdot 5! (16/15) \vol(P_\lambda \cap P_C)$.
In the special cases $\lambda = 0, 1/2, 1$ of plurality, Borda,
antiplurality, respectively, we obtain $119/135$, $41/45$, $17/27$.
These last three results were obtained long ago by Gehrlein.

We can consider further intersections of such events. For example,
Gehrlein has computed limiting results under IC for the
conditional probability that rule $\lambda$ chooses the Condorcet winner
given that Borda does, that Borda does given that rule $\lambda$ does,
and that both rules choose the Condorcet winner given that it exists.
The answers to these questions are easily found for IAC using the above
methods and are listed in Table~\ref{table:cond}. These have not
previously been published as far as we are aware (numbers in brackets in that table 
represent citations). In
Table~\ref{table:cond} we let $A | C$ denote the event that
antiplurality chooses the Condorcet winner given that it exists, $B | (P
\cap C)$ the probability that Borda chooses the Condorcet winner given
that plurality does, etc. These can be computed easily using the events
$C$ and $X_\lambda$ above. For example, the entry $B | (P \cap C)$
corresponds to the probability of the event that Borda and Condorcet
agree given that plurality and Condorcet agree. This is simply the
volume of the polytope $P_{1/2} \cap P_C \cap P_0$ divided by the volume
of $P_0 \cap P_C$ (the factor of $3$ cancels out because we are
computing conditional probabilities via $P(E_1 | E_2) = P(E_1 \cap
E_2)/P(E_1)$).

\begin{table}
\begin{tabular}{|c|c|c|c|c|c|c|c|}
\hline
$P | C$ & $A | C$ & $B | C$ & $(A \cap B) | C$ & $(A \cap P) | C$ & $(B \cap P) | C$ & 
 $B | (P \cap C)$ & $B| (A \cap C)$ \\
\hline
 0.88148 \cite{Gehr1982}& 0.62963 \cite{Gehr1982}& 0.91111\cite{Gehr1992} & 0.61775 & 0.53040 & 0.81821 &  0.92282 & 0.98113 \\
\hline
\end{tabular}
\caption{(Joint) limiting Condorcet efficiencies of the standard positional rules 
under IAC, 3 candidates}
\label{table:cond}
\end{table}

We consider even more intersections of such events in the next section.

In \cite{CGZ2005} the value of $\lambda$ for which the
positional rule with weights $(1, \lambda, 0)$ is most Condorcet
efficient was determined. We call this ``rule M" for brevity.
The optimal value of $\lambda$ is an algebraic irrational number given as
the root of a polynomial of degree $8$ and to $5$ decimal places equals
$0.37228$. The corresponding value of the Condorcet efficiency is
approximately $0.92546$, only slightly more than that for Borda. 
To use this particular value of $\lambda$ in computations similar to
those above, it is probably best to switch to software that performs
floating point computations in order to compute volumes. One such is
{\tt vinci}. We obtain for example that the joint Condorcet efficiency of
the optimal rule and the Borda rule equals, to 5 decimal places, $0.89183$.

\subsubsection*{Borda's Paradox}

We finish here by discussing \Em{Borda's Paradox}. Some rules can elect
a Condorcet loser, namely a candidate that is beaten by every other when
pairwise comparisons are made. The probability of this event for
plurality and antiplurality has been studied under IAC in
\cite{Lepe1993}, and it has long been known to be zero for
Borda. The methods in this section can be applied directly, since we
need only replace the Condorcet winner conditions by the same ones with
the direction of the inequality reversed. This shows that Borda's
Paradox occurs for plurality with probability $1/36$, agreeing with
\cite{Lepe1993}. The corresponding results for Borda and
antiplurality are $0$ and $17/576$, corroborating the previous results.
The probability that the most Condorcet efficient rule above elects the
Condorcet loser is, as one might expect, very small. The results are shown in 
Table~\ref{table:loser}.

\begin{table}
\begin{tabular}{|c|c|c|c|}
\hline
plurality & rule M & Borda & antiplurality \\
\hline
0.0278 \cite{Lepe1993} & 0.00131 & 0 & 0.0295 \cite{Lepe1993}\\
\hline
\end{tabular}
\caption{Limiting probability of Borda's paradox under IAC, 3 candidates}
\label{table:loser}
\end{table}

\subsection{When do all common rules elect the same winner?}
\label{ss:allsame}

For three-alternative elections, all positional voting rules elect the
same winner in a given situation if and only if both plurality and
antiplurality elect the same winner, since the vector of scores is a
convex combination of those for the two extreme rules. The probability
of this event has been investigated under IC but not under IAC as far as
we are aware. In \cite{MTV2000} Merlin, Tataru and
Valognes also investigated the probability under IC that all positional
rules and all Condorcet efficient rules yield the same winner (in this
case, all scoring runoff rules also yield this same winner).

We again suppose that $a$ is the winner. We want to compute the
probability of the event $P \cap A$ as described in the previous
section. The relevant polytope has $18$ vertices and $m = 12$. Its
volume is $113/77760$ and so the limiting probability that all
positional rules yield the same winner for $3$ alternatives under IAC is
$113/216$ (this confirms a result in \cite{Gehr2002b}).
We could also investigate the relationship between, say, plurality and
Borda. They agree with probability $89/108$, whereas
antiplurality and Borda agree with probability $1039/1512$.

The probability that all Condorcet rules and all positional rules elect
the same winner given that the Condorcet winner exists is obtained
easily via computation of $\Pr(P \cap C \cap A)$ as above. The answer is
$3437/6480$. The polytope involved has $29$ vertices and $m = 12$.

We must also consider the case when no Condorcet winner exists. There
are two cases corresponding to the two cycles $a, b, c, a$ and $a, c, b,
a$. In the first case, \cite{MTV2000} shows that the
rules all agree if and only if all positional rules give the ranking $a,
b, c$, and this occurs if and only if both plurality and antiplurality
give that ordering. The computation is straightforward as above and the
probability of this event is only $5/10368$. The contribution from the
cyclic case is therefore 32 times this, or, $5/324$, and the final
result for the probability that all rules agree is $10631/20736$.

We can also consider the probability that two rules agree in their whole
ranking, not just in the choice of winner. This is easily computed
similarly to above: plurality and antiplurality agree on their whole
ranking with probability $8/27$, while Borda and
plurality agree with probability $61/108$.  Borda and
antiplurality also agree with probability $61/108$, which is clear by
symmetry in any case.

\begin{table}
\begin{tabular}{|c|c|c|}
\hline
Rules & Elect same winner & Agree whole ranking\\
\hline
Antiplurality and Borda & 0.68717 &  0.56481 \\
\hline
Antiplurality and plurality (hence all scoring rules) & 0.52315 \cite{Gehr2002b} & 0.29630\\
\hline
Plurality and Borda & 0.82407 &  0.56481\\
\hline
All common rules & 0.51268 & \\
\hline
\end{tabular}
\caption{Limiting probability of agreement of various rules under IAC, 3 candidates}
\label{table:agree}
\end{table}

\subsection{Abstention and Participation Paradoxes}

In \cite{LeMe2001} Lepelley and Merlin discuss various
ways in which voters can attempt to manipulate an election by abstaining
from voting. All scoring runoff rules and Condorcet rules suffer from
this problem. Although abstaining turns out to be a dominated strategy
for scoring runoff rules, it is still of interest to compute the
probability that a situation may be manipulated in this way. Lepelley
and Merlin carry this out under IC and IAC for scoring runoff rules
based on plurality, antiplurality and Borda. For the latter (the Nanson
rule) the limiting probability was not computed exactly (Table~5 of
\cite{LeMe2001} refers to results of Monte Carlo
simulation). We compute some exact values here.

We use the linear system given in \cite{LeMe2001}. Suppose
that $c$ is eliminated first and $a$ then beats $b$ in the runoff. The
Positive Participation Paradox occurs when voters ranking $a$ first are
added to the electorate, and yet $a$ then loses. This cannot happen when
plurality is used at the first stage, but for other rules it can happen
that $b$ now loses the first stage, and $a$ subsequently loses the
runoff against $c$. Note that only voters with preference order $acb$
can cause this to occur, and it can only occur when $c$ originally beats
$a$ pairwise.

The system describing this set of voting situations contains the
inequalities stating that $a$ beats $b$ and $b$ beats $c$ using the
given scoring rule, and also that $a$ beats $b$ pairwise. In addition we
have another constraint as described in \cite{LeMe2001}
(note that $n_6$ in the first equation on p.58 of that paper should be
$-n_6$). Carrying out the (by now routine) computation we obtain 
$1/72$ which confirms the simulation result $0.14$ referred to
above. Note that the polytope involved has only $6$ vertices and $6$
facets but $m =18$; if $e = 18$ (which we have not checked),  it would
be difficult to compute the Ehrhart polynomial using the old methods,
which probably explains why only simulation results were obtained for the Nanson 
rule in the paper cited above.

Similarly we may compute the result for each of several other participation paradoxes. The results
for the negative participation, positive abstention and negative
abstention paradoxes (see \cite{LeMe2001} for definitions
and characterizations of the polytopes) are respectively $1/48, 1/96,
1/72$ confirming the earlier simulation results $0.020, 0.010, 0.14$.

We can also perform the analogous computations for plurality and antiplurality runoff --- the 
results confirm those in \cite{LeMe2001} and are shown in Table~\ref{table:abs}.

\begin{table}
\begin{tabular}{|c|c|c|c|c|}
\hline
Underlying rule & PPP & NPP & PAP & NAP\\
\hline
plurality \cite{LeMe2001} & 0 &0.07292 & 0 & 0.04080\\
Borda & 0.01389 & 0.02083 & 0.01042 & 0.01389 \\
antiplurality \cite{LeMe2001} & 0.03822 &0  &0.04253 & 0\\
\hline
\end{tabular}
\caption{Limiting probability of participation paradoxes for scoring runoff rules under IAC, 
3 candidates}
\label{table:abs}
\end{table}

\subsection{The referendum paradox}
\label{ss:ref}

This gives an example where the variables describing our polytopes are 
slightly different.

In \cite{FLM2004} the referendum or Compound Majority Paradox
is studied.  In the simplest case there are $N$ equal sized districts
each having $n$ voters.  There are two candidates $a$ and $b$ and voters
in each district use majority rule to decide which candidate wins each
district. The candidate winning a majority of districts is the winner of
the election; the paradox occurs when this candidate would have lost if
simple majority had been used in the union of all districts.

Among other things, the authors of \cite{FLM2004} derive the
probability of occurrence for $N=3, 4, 5$ under IAC using the older
methods and state that they are not able to extend it to $N\geq 6$.
Using the methods of the present paper it is easy to perform the
computations for at least a few more values of $N$. Let $n_i$ denote the
number of voters voting for $a$ in district $i$. The relevant set turns
out to be described (ignoring ties for simplicity) by the union of
polytopes of the form
\begin{align*}
n_i & \geq N/2 \text{ for $1 \leq i \leq k$} 
\qquad\text{($a$ wins $k$ districts)} \\
0  \leq n_i & \leq N/2 \text{ for  $k+1 \leq i \leq N$} 
\qquad\text{($b$ wins $N-k$ districts)} \\
\sum_i n_i & \leq Nn/2 \qquad\text{($b$ wins overall)}
\end{align*}
for $\lfloor N/2 \rfloor + 1 \leq k \leq N - 1$. The polytope $P_n$
corresponding to $k$ must have volume multiplied by $2 \binom{N}{k}$ to
account for the symmetries of the problem.  Note that there are $(n+1)^N$
situations to consider and so the leading term in $n$ of the Ehrhart
polynomial gives the probability required.

Doing the analogous computation  for $N = 3$ and $N = 4$ we obtain  (as
$n \to \infty$, in other words computing the volume of $P_n/n$) results
agreeing with \cite{FLM2004}. However already for $N = 5$ we
obtain $61/384$ as opposed to their result $55/384$. For $N = 7$ we have
$9409/46080$. The most complicated corresponding polytope in the last
case has $36$ vertices and $11$ facets, whereas for $N = 9$ it has $91$ vertices and $14$ facets. 
We did not attempt to find the
maximum value of $N$ for which our software could obtain an answer; 
the answer was  given essentially instantaneously for $N = 9$. The conjecture in \cite{FLM2004} that 
the probability tends to a limit of around $0.165$ as $N$ (odd) goes to infinity seems unlikely 
in the light of these results.

\begin{table}
\begin{tabular}{|c|c|c|c|c|c|c|}
\hline
number of districts & 3 & 4 & 5 & 6 & 7 & 9\\
\hline
probability &0.125 \cite{FLM2004} &0.02083 \cite{FLM2004} &0.15885 \cite{FLM2004} 
&0.04063 & 0.20419 & 0.26954 \\
\hline
\end{tabular}
\caption{Limiting probability of referendum paradox under IAC, 
3 candidates}
\label{table:ref}
\end{table}

\section{Summary and discussion of future work}
We have shown that a wide variety of natural probabilistic questions for
$3$-alternative elections under IAC can be answered by applying standard
algorithms for counting lattice points in, and computing volumes of,
convex polytopes. For $4$ or more alternatives the computations are
conceptually the same but necessarily more complicated. However, the
scope for extending results in the $3$-candidate case to $4$ or more
candidates is obviously higher than for the older methods, which now
appear to be completely superseded. One important point to notice is
that many algorithms for volume computation have running times that are
very sensitive to the number of defining hyperplanes and the number of
vertices. Thus finding the most efficient description of the input
system is important. It is certainly clear that further progress in this
area will require researchers in social choice theory to understand in
some detail how the fastest algorithms for lattice point counting and
volume computation actually work. This may even lead to proofs for
larger (or general) numbers of candidates when the polytopes concerned
have a particularly nice structure.

Many questions naturally arise from our work here. One obvious line of
attack is to try to find the optimal parameter for $3$-alternative scoring
rules that minimizes the probability of a certain undesirable behaviour
occurring. The present authors are already engaged in carrying this out
for the case of (naive, coalitional) manipulability. Numerical results
obtained in \cite{PrWi2005} show that the answer may well be
plurality, but this has never been proved. An attack on this problem
along the lines of the approach in the present paper would require
computation of volumes of a polytope whose defining constraints depend
linearly in a parameter $\lambda$, and this requires considerable work
as shown in \cite{CGZ2005}. Understanding of how to
carry out such a computation would help in understanding the variation
between positional rules. For example, the probability of electing a
Condorcet loser is of order $0.03$ for both plurality and antiplurality,
but an order of magnitude smaller near the Borda rule, and as a function
of $\lambda$ is very flat there. Quantifying this type of variation
analytically may show, for example, that it is not worth the trouble of
replacing Borda by the Condorcet-optimal positional rule.

Another direction is to consider other probability models. For
simplicity here we have not considered some common assumptions such as
single-peaked preferences and the Maximal Culture Condition. Many
computations in these cases reduce to ones identical in spirit to those
we have undertaken here. More general P{\'o}lya-Eggenberger distributions
would lead to the more difficult issue of integrals of nonconstant
probability densities over polytopes, but some results may be
forthcoming there.

\bibliographystyle{plain}
\bibliography{voting}

\end{document}